\documentclass{article}
\usepackage[letterpaper,top=3.5cm,bottom=3.5cm,left=2.5cm,right=2.5cm,marginparwidth=2.75cm]{geometry}

\usepackage{mathtools}
\usepackage{amsmath}
\usepackage{amssymb}
\usepackage{amsthm}
\usepackage[square,numbers]{natbib}
\usepackage{xparse}

\newtheorem{theorem}{Theorem}[section]

\newtheorem{corollary}[theorem]{Corollary}

\newtheorem{lemma}[theorem]{Lemma}

\newtheorem{proposition}[theorem]{Proposition}

\DeclareMathOperator{\re}{\mathfrak{Re}}
\DeclareMathOperator{\res}{Res}
\makeatletter
\NewDocumentCommand{\sump}{e{_}}
 {%
  \DOTSB
  \mathop{\IfNoValueTF{#1}{\sump@{}}{\sump@{#1}}}%
  \nolimits
 }
\newcommand{\sump@}[1]{\mathpalette\sump@@{#1}}
\newcommand{\sump@@}[2]{%
  \ifx#1\displaystyle
    {\sump@display{#2}}%
  \else
    \sum@\nolimits'_{#2}%
  \fi
}
\newcommand{\sump@display}[1]{%
  \sbox\z@{$\m@th\displaystyle\sum@\nolimits'$}%
  \sbox\tw@{$\m@th\displaystyle\sum@\limits_{#1}$}%
  \sbox\@tempboxa{$\m@th\displaystyle'$}
  \mathop{\sum@\nolimits' \kern-\wd\@tempboxa}\limits_{#1}%
  \ifdim\wd\z@>\wd\tw@
    \kern\dimexpr\wd\z@-\wd\tw@\relax
  \fi
}
\makeatother

\title{A Product Identity For Dirichlet Series Satisfying Hecke's Functional Equation}
\author{Efe Gürel$^{1}$\\
        \small $^{1}$TÜBİTAK Natural Sciences High School, Kocaeli, 41400, Turkey}
\date{}

\begin{document}

\maketitle
\begin{abstract}
    In this paper, we give an analogue of Wilton's product formula for Dirichlet series that satisfy Hecke's functional equation. We apply our results to obtain identities for Hecke series, L-functions associated to modular forms, Ramanujan's L-function, Epstein zeta functions, Dedekind zeta functions of imaginary quadratic fields and Dirichlet L-functions. A $4$-term product identity for Riemann zeta function is also given.\\

     

\end{abstract}

\section{Introduction}

The Riemann zeta function is defined as
\begin{align*}
    \zeta(s)=\sum_{n=1}^{\infty}\frac{1}{n^s} \qquad (\re(s)=\sigma>1)
\end{align*}
and has an analytic continuation to entire complex plane except $s=1$ where it has a simple pole. It satisfies the remarkable functional equation
\begin{align}\label{RiemannFe}
    \pi^{-s/2}\Gamma\Big(\frac{s}{2}\Big)\zeta(s)=\pi^{(1-s)/2}\Gamma\Big(\frac{1-s}{2}\Big)\zeta(1-s)
\end{align}
due to Riemann. Here $\Gamma$ denotes the well-known gamma function. Approximate functional equations are fundamental objects for studying the behavior of zeta functions near the critical strip, such as obtaining mean value theorems. Approximate functional equations were obtained by Hardy and Littlewood in \cite{HardyApprox} for $\zeta(s)$ and $\zeta(s)^2$. Wilton \cite{wiltonOG} investigated the approximate functional for the product of two zeta functions. His method relies on the following theorem which is now known as Wilton's formula \cite{Linearized}.
\begin{theorem}
    Let $u,v\in\mathbb{C}$ such that $\re(u),\re(v)>-1$, $\re(u+v)>0$, $u,v\neq 1$ and $u+v\neq 2$. Then we have
    \begin{align*}
        \zeta(u)\zeta(v)-\left( \frac{1}{u-1}+\frac{1}{v-1} \right)\zeta(u+v-1)&=2(2\pi)^{u-1}\sum_{n=1}^{\infty}\sigma_{1-u-v}(n)n^{u-1}u\int_{2\pi n}^{\infty}t^{-u-1}\sin t dt\\
        &+2(2\pi)^{v-1}\sum_{n=1}^{\infty}\sigma_{1-u-v}(n)n^{v-1}v\int_{2\pi n}^{\infty}t^{-v-1}\sin t dt,
    \end{align*}
    where $\sigma_z(n)=\sum_{d|n}d^z$.
\end{theorem}
An analogue of Wilton’s formula for Dedekind zeta functions and applications were discussed in \cite{DedekindWilton2,DedekindWilton3,DedekindWilton} and references therein. 
\newline

An analytic function defined on the complex upper half-plane $\mathbb{H}$ is said to be a modular form of order $k\in \mathbb{Z}$ if it has a Fourier expansion of form $f(\tau)=\sum_{n=0}^\infty a_n e^{2\pi i n\tau}$ and satisfies $f\left( \frac{a\tau+b}{c\tau+d} \right)=(c\tau+d)^k f(\tau)$ whenever $\left(\begin{smallmatrix}a & b \\c & d\end{smallmatrix}\right)\in SL_2(\mathbb{Z})$ \cite{ApostolModular}. L-functions associated with modular forms are of great importance. Hecke's correspondence theorem \cite{BerndtHeckeTheory} connects functions satisfying a modular relation with L-functions' functional equations as follows. 
\begin{theorem}\label{HCT}
    Let $\lambda>0$, $k\in\mathbb{R}$ and $\gamma\in\mathbb{C}$. Furthermore, let $\left\{ \alpha_n \right\}$ and $\left\{ \beta_n \right\}$ be sequences of complex numbers such that $\alpha_n,\beta_n=O(n^c)$ for some $c\ge 0$. Define for $\sigma>c+1$
    \begin{align*}
        \varphi(s)=\sum_{n=1}^{\infty}\alpha_n n^{-s} \qquad \text{and}\qquad \psi(s)=\sum_{n=1}^{\infty}\beta_n n^{-s}
    \end{align*}
    and
    \begin{align*}
        \Phi(s)=\left( \frac{2\pi}{\lambda} \right)^{-s}\Gamma(s)\varphi(s) \qquad \text{and}\qquad \Psi(s)=\left( \frac{2\pi}{\lambda} \right)^{-s}\Gamma(s)\psi(s).
    \end{align*}
    For $\tau\in\mathbb{H}$, let
    \begin{align*}
        f_\alpha(\tau)=\sum_{n=0}^{\infty}\alpha_n e^{2\pi i n\tau/\lambda} \qquad \text{and}\qquad f_\beta(\tau)=\sum_{n=0}^{\infty}\beta_n e^{2\pi i n\tau/\lambda}.
    \end{align*}
    Then the following two statements are equivalent.
    \begin{enumerate}
        \item $f_\alpha(\tau)=\gamma\left( \tau/i \right)^{-k}f_\beta(-1/\tau).$
        \item The functions
        \begin{align*}
            \Phi(s)+\frac{\alpha_0}{s}+\frac{\gamma\beta_0}{k-s} \qquad \text{and}\qquad \Psi(s)+\frac{\beta_0}{s}+\frac{\alpha_0}{\gamma(k-s)}
        \end{align*}
        have an analytic continuation to the entire complex plane that are entire and bounded in every vertical strip. Moreover,
        \begin{align*}
            \Phi(s)=\gamma\Psi(k-s).
        \end{align*}
    \end{enumerate}
\end{theorem}
If $\psi=\varphi$ and $\varphi$ satisfies one of the conditions in Theorem \ref{HCT}, then $\varphi$ is said to be a Hecke series of Hecke signature $(\lambda,k,\gamma)$. A Hecke series not identical to zero has an automorphy factor of $\gamma=\pm 1$. When $k\ge4$ is an even natural number, Hecke series of signature $(1,k,(-1)^{k/2})$ correspond to the modular forms of weight $k$. 
\newline

Approximate functional equations for Hecke series were investigated by Apostol in \cite{ApostolHecke} and for more general classes of Dirichlet series in \cite{ChandraApprox,LavrikApprox}. 
\newline

The Meijer G function \cite{HigherTranscendental} is defined as the Mellin-Barnes integral 
\begin{align} \label{MeijerDef}
    G_{p,q}^{m,n}\left( \begin{array}{c} a_1,\ldots,a_p\\ b_1,\ldots,b_q \end{array}\middle\vert z \right)=\frac{1}{2\pi i}\int_{L}\dfrac{\prod_{j=1}^{m}\Gamma(b_j-s)\prod_{j=1}^{n}\Gamma(1-a_j+s)}{\prod_{j=m+1}^{q}\Gamma(1-b_j+s)\prod_{j=n+1}^{p}\Gamma(a_j-s)} z^s ds,
\end{align}
where $m,n,p,q$ are non-negative integers such that $m\le q$, $n\le p$ with $z,a_1,\ldots a_p,b_1,\ldots b_q\in\mathbb{C}$ and the poles of the functions $\Gamma(b_j-s)$ and $\Gamma(1-a_k+s)$ do not coincide for any $j=1,\ldots,m$ and $k=1,\ldots,n$. Here $L$ is a line that goes from $-i\infty$ to $i\infty$. For the integral in \eqref{MeijerDef} to converge it suffices if $p+q<2(m+n)$ and $|\arg z|<\pi\left(m+n-\frac{p+q}{2}\right)$. The integral also converges for $|\arg z|=\pi\left(m+n-\frac{p+q}{2}\right)$ if
\begin{align*}
    (q-p)\left( \re (s)+\frac{1}{2} \right)>1+\sum_{j=1}^{q}\re (b_j) -\sum_{j=1}^{p}\re (a_j).
\end{align*}
One special case we need is 
\begin{align*}
    G_{0,2}^{1,0}\left( \begin{array}{c},\\ a,b \end{array}\middle\vert z \right)= z^{\frac{a+b}{2}}J_{a-b}\left( 2z^{1/2} \right)
\end{align*}
where $J_\rho(z)$ is the Bessel function of first kind and order $\rho$, given by
\begin{align*}
    J_\rho(z)=\sum_{n=0}^{\infty}\frac{(-1)^n}{\Gamma(n+1)\Gamma(n+\rho+1)}\left( \frac{z}{2} \right)^{2n+\rho}.
\end{align*}
\newline

We state our main result as follows.
\begin{theorem}\label{MainThm}
    Let $\varphi(s),\psi(s), (\lambda,k,\gamma)$ and $c$ satisfy the conditions and one of the equivalent statements of Theorem \ref{HCT}. Then, for $\re(u),\re(v)>\max(c+1,k)$ and $u,v\neq k+1$, we have
    \begin{align}\label{MainThmEq}
        \begin{split}
            \varphi(u)\psi(v)&=\frac{\res\varphi(k)}{u-k}\psi(u+v-k)+\frac{\res\psi(k)}{v-k}\varphi(u+v-k)\\
            &-\frac{2\pi\gamma}{\lambda}\sum_{n=1}^{\infty}\sigma_{\beta,k-u-v}(n) n^{\frac{1-k}{2}}\int_{0}^{1}t^{\frac{k-1}{2}-u}J_{k-1}\left( \frac{4\pi\sqrt{nt}}{\lambda} \right)dt\\
            &-\frac{2\pi}{\lambda\gamma}\sum_{n=1}^{\infty}\sigma_{\alpha,k-u-v}(n) n^{\frac{1-k}{2}}\int_{0}^{1}t^{\frac{k-1}{2}-v}J_{k-1}\left( \frac{4\pi\sqrt{nt}}{\lambda} \right)dt
        \end{split} 
    \end{align}
    where
    \begin{align*}
        \sigma_{\alpha,z}(n)=\sum_{d|n}^{}\alpha_d \alpha_{n/d} d^z \qquad \text{and}\qquad \sigma_{\beta,z}(n)=\sum_{d|n}^{}\beta_d \beta_{n/d} d^z.
    \end{align*}

\end{theorem}
The residues of $\varphi$ and $\psi$ at $s=k$ can be easily computed by Theorem \ref{HCT} and are given by
\begin{align*}
        \res\varphi(k)=\left( \frac{2\pi}{\lambda} \right)^k \frac{\gamma\beta_0}{\Gamma(k)} \qquad \text{and}\qquad \res\psi(k)=\left( \frac{2\pi}{\lambda} \right)^k \frac{\alpha_0}{\gamma\Gamma(k)}.
    \end{align*}
We remark that by analytic continuation, equation \eqref{MainThmEq} may be valid for larger domains of $u$ and $v$. We also note that equation \eqref{MainThmEq} is highly suitable for computing special values of $\varphi$ and $\psi$.
\section{Preliminaries}

Our method is based on a special case of Riesz summation, a generalization of Perron's formula, and is similar to that of \cite{Linearized,DedekindWilton}. We refer the reader to \cite{Riesz} for a discussion of the general theory of Dirichlet series and Riesz summation. Let $\varphi(s)=\sum_{n=1}^{\infty}\alpha_n n^{-s}$ for $\sigma>\sigma_\varphi$ and $ \psi(s)=\sum_{n=1}^{\infty}\beta_n n^{-s}$ for $\sigma>\sigma_\psi$ be Dirichlet series that can be meromorphically continued to the complex plane and satisfy suitable decaying conditions.  Define the integral operator
\begin{align*}
    \mathcal{F}_{a}\left( \varphi(u),\psi(v);x \right)=\frac{1}{2\pi i}\int_{a-i\infty}^{a+i\infty}\varphi(u+w)\psi(v-w)\frac{x^{w+1}}{w(w+1)}dw,
\end{align*}
where $\re(u)>\sigma_\varphi+a$, $\re(v)>\sigma_\psi+a$ and $x$ is a real variable. This operator corresponds to a special case of the general Perron's formula. The conditions of meromorphicity and growth are assured by Theorem \ref{HCT} and $\sigma_\varphi=\sigma_\psi=c+1$ in our case. We now use Perron's formula to obtain
\begin{align}\label{FPhiPsiSeries}
\begin{split}
    \mathcal{F}_{a}\left( \varphi(u),\psi(v);x \right)&=\frac{1}{2\pi i}\int_{a-i\infty}^{a+i\infty}\varphi(u+w)\psi(v-w)\frac{x^{w+1}}{w(w+1)}dw\\
&=\frac{1}{2\pi i}\int_{a-i\infty}^{a+i\infty}\frac{x^{w+1}}{w(w+1)}\sum_{n=1}^{\infty}\frac{\alpha_n}{n^{u+w}}\sum_{m=1}^{\infty}\frac{\beta_m}{m^{v-w}}dw\\
&=\sum_{m=1}^{\infty}\frac{\beta_m}{m^v}\sum_{n=1}^{\infty}\frac{\alpha_n}{n^u}\frac{x}{2\pi i}\int_{a-i\infty}^{a+i\infty}\left( \frac{n}{mx} \right)^{-w}\frac{dw}{w(w+1)}\\
&=\sum_{m=1}^{\infty}\frac{\beta_m}{m^{v+1}}\sump_{n\le mx} \frac{\alpha_n}{n^u}\left( mx-n \right),
\end{split}
\end{align}
where $\sump_{\ell\le x}$ indicates that the last term is to be halved when $x$ is an integer. Here the interchange of integration and summations are justified by the absolute convergence. In a similar manner, we have
\begin{align}\label{FPsiPhiSeries}
    \mathcal{F}_{a}\left( \psi(v),\varphi(u);x \right)=\sum_{n=1}^{\infty}\frac{\alpha_n}{n^{u+1}}\sump_{m\le nx} \frac{\beta_m}{m^v}\left( nx-m \right).
\end{align}
Our main ingredient is a Perron-like formula for $\varphi(u)\psi(v)$, which has been derived for special cases in \cite{Linearized,DedekindWilton}. 
\begin{lemma}\label{PerronLemma}
    Let $a>0$, $\re(u)>\sigma_\varphi+a$ and $\re(v)>\sigma_\psi+a$. Then, the following holds:
    \begin{align*}
        \frac{\partial }{\partial x}\mathcal{F}_{a}\left( \varphi(u),\psi(v);1 \right)+\frac{\partial }{\partial x}\mathcal{F}_{a}\left( \psi(v),\varphi(u);1 \right)=\varphi(u)\psi(v).
    \end{align*}
\end{lemma}
\begin{proof}
    Differentiating equation \eqref{FPhiPsiSeries} with respect to $x$ and putting $x=1$, we obtain
    \begin{align*}
         \frac{\partial }{\partial x}\mathcal{F}_{a}\left( \varphi(u),\psi(v);1 \right)=\sum_{m=1}^{\infty}\sump_{n\le m} \frac{\beta_m}{m^{v}}\frac{\alpha_n}{n^u}
    \end{align*}
    and similarly
    \begin{align*}
         \frac{\partial }{\partial x}\mathcal{F}_{a}\left( \psi(v),\varphi(u);1 \right)=\sum_{n=1}^{\infty}\sump_{m\le n} \frac{\beta_m}{m^{v}}\frac{\alpha_n}{n^u}.
    \end{align*}
    Adding these together we have
    \begin{align*}
         \frac{\partial }{\partial x}\mathcal{F}_{a}\left( \varphi(u),\psi(v);1 \right)+\frac{\partial }{\partial x}\mathcal{F}_{a}\left( \psi(v),\varphi(u);1 \right)&=\sum_{m=1}^{\infty}\sump_{n\le m} \frac{\beta_m}{m^{v}}\frac{\alpha_n}{n^u}+\sum_{n=1}^{\infty}\sump_{m\le n} \frac{\beta_m}{m^{v}}\frac{\alpha_n}{n^u}\\
        &=\sum_{m=1}^{\infty}\sum_{n=1}^{\infty}\frac{\beta_m}{m^{v}}\frac{\alpha_n}{n^u}\\
        &=\varphi(u)\psi(v),
    \end{align*}
    where we have used the so called Nakajima disection $
\sum_{m=1}^{\infty}\sum_{n=1}^{\infty}=\sum_{m=1}^{\infty}\sump_{n\le m}+\sum_{n=1}^{\infty}\sump_{m\le n}$ \cite{Nakajima}. This completes the proof.
\end{proof}

We need the following lemma.
\begin{lemma}\label{DirichletMultipLemma}
    For $\re(z)>c+1-\re(u), c+1+\re(v)-k$, we have
    \begin{align*}
        \varphi(u+z)\varphi(k-v+z)=\sum_{n=1}^{\infty}\frac{\sigma_{\alpha,k-u-v}(n)}{n^{z+k-v}}.
    \end{align*}
    Furthermore, for $\re(z)>c+1-\re(v), c+1+\re(u)-k$, we have
    \begin{align*}
        \psi(v+z)\psi(k-u+z)=\sum_{n=1}^{\infty}\frac{\sigma_{\beta,k-u-v}(n)}{n^{z+k-u}}.
    \end{align*}
\end{lemma}
\begin{proof}
    By the condition on $\re(z)$, both of the series $\varphi(u+z)$ and $\varphi(k-v+z)$ converge absolutely. Thus we get
    \begin{align*}
        \varphi(u+z)\varphi(k-v+z)&=\sum_{d=1}^{\infty}\sum_{m=1}^{\infty}\frac{\alpha_d}{d^{u+z}}\frac{\alpha_m}{m^{k-v+z}}\\
        &=\sum_{n=1}^{\infty}\sum_{dm=n}\frac{\alpha_d}{d^{u+z}}\frac{\alpha_m}{m^{k-v+z}}\\
        &=\sum_{n=1}^{\infty}\frac{1}{n^{z+k-v}}\sum_{d|n}\alpha_d\alpha_{n/d}d^{k-u-v}\\
        &=\sum_{n=1}^{\infty}\frac{\sigma_{\alpha,k-u-v}(n)}{n^{z+k-v}}.
    \end{align*}
    The proof for $\psi$ is completely analogous. This completes the proof.
\end{proof}
\section{Proof of Theorem \ref{MainThm}}
Take some $1\neq a>0$ such that
\begin{align*}
    a>\max(c+1-\re(u),c+1-\re(v), c+1+\re(u)-k,c+1+\re(v)-k).
\end{align*}
By the assumptions of Theorem \ref{MainThm}, we have
\begin{align*}
    \varphi(s)=\gamma\left( \frac{2\pi}{\lambda} \right)^{2s-k}\frac{\Gamma(k-s)}{\Gamma(s)}\psi(k-s).
\end{align*}
Applying the above functional equation in the definition of $\mathcal{F}_{-a}$, we get
\begin{align}\label{F-aWtoZ}
    \begin{split}
        \mathcal{F}_{-a}\left( \varphi(u),\psi(v);x \right)&=\frac{1}{2\pi i}\int_{-a-i\infty}^{-a+i\infty}\varphi(u+w)\psi(v-w)\frac{x^{w+1}}{w(w+1)}dw\\
        &=\frac{\gamma}{2\pi i}\int_{-a-i\infty}^{-a+i\infty}\left( \frac{2\pi}{\lambda} \right)^{2u+2w-k}\frac{\Gamma(k-u-w)}{\Gamma(u+w)}\psi(v-w)\psi(k-u-w)\frac{x^{w+1}}{w(w+1)}dw\\
        &=\frac{\gamma}{2\pi i}\int_{a-i\infty}^{a+i\infty}\left( \frac{2\pi}{\lambda} \right)^{2u-2z-k}\frac{\Gamma(k-u+z)}{\Gamma(u-z)}\psi(v+z)\psi(k-u+z)\frac{x^{1-z}}{z(z-1)}dz,
    \end{split}
   \end{align}
where we changed the variables $z=-w$. Since $a=\re(z)>c+1-\re(v), c+1+\re(u)-k$, we can substitute Lemma \ref{DirichletMultipLemma} for $\psi(v+z)\psi(k-u+z)$ in \eqref{F-aWtoZ} to obtain
\begin{align*}
    \mathcal{F}_{-a}\left( \varphi(u),\psi(v);x \right)=\frac{\gamma}{2\pi i}\int_{a-i\infty}^{a+i\infty}\left( \frac{2\pi}{\lambda} \right)^{2u-2z-k}\frac{\Gamma(k-u+z)}{\Gamma(u-z)}\frac{x^{1-z}}{z(z-1)}\sum_{n=1}^{\infty}\frac{\sigma_{\beta,k-u-v}(n)}{n^{z+k-u}}dz.
\end{align*}
Since $\beta_n=O(n^c),$ we easily get $\sigma_{\beta,s}(n)=O\left(n^c \sigma_{\re(s)}(n)\right)=O\left(n^{c+\re(s)}\right)$. From this estimate, we have that $\sigma_{\beta,k-u-v}(n)/n^{z+k-u}=O\left(n^{c-a-\re(v)}\right)$ which is $O\left(n^{-1-\epsilon}\right)$ for some $\epsilon>0$ by the condition on $a$. Therefore the interchange of orders of operators are justified by the absolute convergence. Interchanging the orders of integral and summation yields
\begin{align*}
    \mathcal{F}_{-a}\left( \varphi(u),\psi(v);x \right)=\gamma\left( \frac{2\pi}{\lambda} \right)^{2u-k}\sum_{n=1}^{\infty}\frac{\sigma_{\beta,k-u-v}(n)}{n^{k-u}}\frac{1}{2\pi i}\int_{a-i\infty}^{a+i\infty}\left( \frac{\lambda^2}{4\pi^2n} \right)^z \frac{\Gamma(k-u+z)}{\Gamma(u-z)}\frac{x^{1-z}}{z(z-1)}dz.
\end{align*}
Differentiating with respect to $x$, we obtain
\begin{align}\label{FirstDxF}
    \frac{\partial }{\partial x}\mathcal{F}_{-a}\left( \varphi(u),\psi(v);x \right)=-\gamma\left( \frac{2\pi}{\lambda} \right)^{2u-k}\sum_{n=1}^{\infty}\frac{\sigma_{\beta,k-u-v}(n)}{n^{k-u}}\frac{1}{2\pi i}\int_{a-i\infty}^{a+i\infty}\left( \frac{\lambda^2}{4\pi^2n} \right)^z \frac{\Gamma(k-u+z)}{\Gamma(u-z)}\frac{x^{-z}}{z}dz.
\end{align}
Now we consider the integrals
\begin{align} \label{I_nDef}
    I_n(x)=\frac{1}{2\pi i}\int_{a-i\infty}^{a+i\infty}\left( \frac{\lambda^2}{4\pi^2n} \right)^z \frac{\Gamma(k-u+z)}{\Gamma(u-z)}\frac{x^{-z}}{z}dz=G_{3,1}^{0,2}\left( \begin{array}{c} 1+u-k,1,u \\ 0 \end{array}\middle\vert \frac{\lambda^2}{4\pi^2 n x} \right).
\end{align}
It is obvious by equation \eqref{FirstDxF} that
\begin{align}\label{SecondDxF}
    \frac{\partial }{\partial x}\mathcal{F}_{-a}\left( \varphi(u),\psi(v);x \right)=-\gamma\left( \frac{2\pi}{\lambda} \right)^{2u-k}\sum_{n=1}^{\infty}\frac{\sigma_{\beta,k-u-v}(n)}{n^{k-u}}I_n(x).
\end{align}
Differentiating equation \eqref{I_nDef} and changing the variables $s=-z$, we have
\begin{align*}
    I'_n(x)&=-\frac{1}{2\pi i x}\int_{-a-i\infty}^{-a+i\infty}\left( \frac{4\pi^2nx}{\lambda^2} \right)^s \frac{\Gamma(k-u-s)}{\Gamma(u+s)}ds\\
    &=-\frac{1}{x}G_{0,2}^{1,0}\left( \begin{array}{c} ,\\ k-u,1-u \end{array}\middle\vert  \frac{4\pi^2nx}{\lambda^2}\right)\\
    &=-\frac{1}{x}\left( \frac{4\pi^2nx}{\lambda^2} \right)^{\frac{k+1}{2}-u}J_{k-1}\left( \frac{4\pi\sqrt{nx}}{\lambda} \right).
\end{align*}
Integrating this result from $0$ to $x$ yields
\begin{align}\label{I_nBessel}
    I_n(x)=-\left( \frac{2\pi}{\lambda} \right)^{k+1-2u}n^{\frac{k+1}{2}-u}\int_{0}^{x}t^{\frac{k-1}{2}-u}J_{k-1}\left( \frac{4\pi\sqrt{nt}}{\lambda} \right)dt.
\end{align}
Combining equations \eqref{SecondDxF} and \eqref{I_nBessel}, we get
\begin{align*}
        \frac{\partial }{\partial x}\mathcal{F}_{-a}\left( \varphi(u),\psi(v);x \right)=\frac{2\pi\gamma}{\lambda}\sum_{n=1}^{\infty}\sigma_{\beta,k-u-v}(n) n^{\frac{1-k}{2}}\int_{0}^{x}t^{\frac{k-1}{2}-u}J_{k-1}\left( \frac{4\pi\sqrt{nt}}{\lambda} \right)dt.
\end{align*}
Similarly, 
\begin{align*}
        \frac{\partial }{\partial x}\mathcal{F}_{-a}\left( \psi(v),\varphi(u);x \right)=\frac{2\pi}{\lambda\gamma}\sum_{n=1}^{\infty}\sigma_{\alpha,k-u-v}(n) n^{\frac{1-k}{2}}\int_{0}^{x}t^{\frac{k-1}{2}-v}J_{k-1}\left( \frac{4\pi\sqrt{nt}}{\lambda} \right)dt.
\end{align*}
Finally this implies
\begin{align}\label{SumOfF-a}
    \begin{split}
        \frac{\partial }{\partial x}\mathcal{F}_{-a}\left( \varphi(u),\psi(v);1 \right)+\frac{\partial }{\partial x}\mathcal{F}_{-a}\left( \psi(v),\varphi(u);1 \right)&=\frac{2\pi\gamma}{\lambda}\sum_{n=1}^{\infty}\sigma_{\beta,k-u-v}(n) n^{\frac{1-k}{2}}\int_{0}^{1}t^{\frac{k-1}{2}-u}J_{k-1}\left( \frac{4\pi\sqrt{nt}}{\lambda} \right)dt\\
        &+\frac{2\pi}{\lambda\gamma}\sum_{n=1}^{\infty}\sigma_{\alpha,k-u-v}(n) n^{\frac{1-k}{2}}\int_{0}^{1}t^{\frac{k-1}{2}-v}J_{k-1}\left( \frac{4\pi\sqrt{nt}}{\lambda} \right)dt.
    \end{split} 
\end{align}
Take some $b>0$ such that
\begin{align*}
    b<\max\left( \re(u)-c-1,\re(v)-c-1,\re(u)-k,\re(v)-k \right).
\end{align*}
Consider the square $S$ with vertices $-a-iT,b-it,b+iT,-a+it$, where $T$ is a sufficiently large real number. Let the function $f$ be defined as
\begin{align*}
    f(w;x)=\varphi(u+w)\psi(v-w)\frac{x^{w+1}}{w(w+1)}.
\end{align*}
In the interior of $S$, $f(w;x)$ has poles at the points $w=0,k-u$ and possibly at $w=-1$. We now integrate $f(w;x)$ on the square $S$. Using Cauchy's residue theorem and taking $T\to \infty$, we obtain
\begin{align*}
    \mathcal{F}_{b}\left( \varphi(u),\psi(v);x \right)&=\mathcal{F}_{-a}\left( \varphi(u),\psi(v);x \right)+x\varphi(u)\psi(v)-\delta_a\varphi(u-1)\psi(v+1)\\
    &+\frac{x^{1+k-u}}{(k-u)(1+k-u)}\res\varphi(k)\psi(u+v-k),
\end{align*}
here $\delta_a$ equals $1$ if $a>1$ and $0$ otherwise. Differentiating with respect to $x$ and putting $x=1$ yields
\begin{align*}
    \frac{\partial }{\partial x}\mathcal{F}_{b}\left( \varphi(u),\psi(v);1 \right)=\frac{\partial }{\partial x}\mathcal{F}_{-a}\left( \varphi(u),\psi(v);1 \right)+\varphi(u)\psi(v)+\frac{\res\varphi(k)}{k-u}\psi(u+v-k).
\end{align*}
Similarly,
\begin{align*}
    \frac{\partial }{\partial x}\mathcal{F}_{b}\left(\psi(v),\varphi(u);1 \right)=\frac{\partial }{\partial x}\mathcal{F}_{-a}\left( \psi(v),\varphi(u);1 \right)+\varphi(u)\psi(v)+\frac{\res\psi(k)}{k-v}\varphi(u+v-k).
\end{align*}
Adding these equations together and using Lemma \ref{PerronLemma}, we have
\begin{align}\label{PhiPsiAlmost}
    \begin{split}
        \varphi(u)\psi(v)&=\frac{\res\varphi(k)}{u-k}\psi(u+v-k)+\frac{\res\psi(k)}{v-k}\varphi(u+v-k)\\
        &-\left( \frac{\partial }{\partial x}\mathcal{F}_{-a}\left( \varphi(u),\psi(v);1 \right)+\frac{\partial }{\partial x}\mathcal{F}_{-a}\left( \psi(v),\varphi(u);1 \right) \right).
    \end{split}
\end{align}
Now substituting equation \eqref{SumOfF-a} into equation \eqref{PhiPsiAlmost}, we get
\begin{align*}
    \varphi(u)\psi(v)&=\frac{\res\varphi(k)}{u-k}\psi(u+v-k)+\frac{\res\psi(k)}{v-k}\varphi(u+v-k)\\
    &-\frac{2\pi\gamma}{\lambda}\sum_{n=1}^{\infty}\sigma_{\beta,k-u-v}(n) n^{\frac{1-k}{2}}\int_{0}^{1}t^{\frac{k-1}{2}-u}J_{k-1}\left( \frac{4\pi\sqrt{nt}}{\lambda} \right)dt\\
    &-\frac{2\pi}{\lambda\gamma}\sum_{n=1}^{\infty}\sigma_{\alpha,k-u-v}(n) n^{\frac{1-k}{2}}\int_{0}^{1}t^{\frac{k-1}{2}-v}J_{k-1}\left( \frac{4\pi\sqrt{nt}}{\lambda} \right)dt.
\end{align*}
This completes the proof. \qed

\section{Applications To Various L-Functions}
In this section we give a few corollaries of the Theorem \ref{MainThm}. Namely, we discuss applications to Hecke series, L-functions associated with modular forms, Ramanujan's L-function, Epstein zeta functions, Dedekind zeta functions of imaginary quadratic fields, Dirichlet L-functions and the Riemann zeta function. Our first obvious corollary is obtained by taking $u=v$.
\begin{corollary}
    For $\re(u)>\max(c+1,k)$, and $u\neq k+1$ we have
    \begin{align*}
            \varphi(u)\psi(u)&=\frac{\res\varphi(k)}{u-k}\psi(2u-k)+\frac{\res\psi(k)}{u-k}\varphi(2u-k)\\
            &-\frac{2\pi\gamma}{\lambda}\sum_{n=1}^{\infty}\sigma_{\beta,k-2u}(n) n^{\frac{1-k}{2}}\int_{0}^{1}t^{\frac{k-1}{2}-u}J_{k-1}\left( \frac{4\pi\sqrt{nt}}{\lambda} \right)dt\\
            &-\frac{2\pi}{\lambda\gamma}\sum_{n=1}^{\infty}\sigma_{\alpha,k-2u}(n) n^{\frac{1-k}{2}}\int_{0}^{1}t^{\frac{k-1}{2}-u}J_{k-1}\left( \frac{4\pi\sqrt{nt}}{\lambda} \right)dt.
    \end{align*}
\end{corollary}
\subsection{Hecke Series}
Let $\varphi$ be a Hecke series of signature $(\lambda,k,\gamma)$ where $\gamma=\pm 1$. Then we have the following.

\begin{proposition}
    For $\re(u),\re(v)>\max(c+1,k)$, and $u,v\neq k+1$, we have
    \begin{align*}\label{MainThmEq}
        \begin{split}
            \varphi(u)\varphi(v)&=\res\varphi(k)\left( \frac{1}{u-k}+\frac{1}{v-k} \right)\varphi(u+v-k)\\
            &-\frac{2\pi\gamma}{\lambda}\sum_{n=1}^{\infty}\sigma_{\alpha,k-u-v}(n) n^{\frac{1-k}{2}}\int_{0}^{1}t^{\frac{k-1}{2}-u}J_{k-1}\left( \frac{4\pi\sqrt{nt}}{\lambda} \right)dt\\
            &-\frac{2\pi\gamma}{\lambda}\sum_{n=1}^{\infty}\sigma_{\alpha,k-u-v}(n) n^{\frac{1-k}{2}}\int_{0}^{1}t^{\frac{k-1}{2}-v}J_{k-1}\left( \frac{4\pi\sqrt{nt}}{\lambda} \right)dt.
        \end{split} 
    \end{align*}
    Furthermore, for $\re(u)>\max(c+1,k)$, and $u\neq k+1$, we have
     \begin{align*}
         \varphi(u)^2=\frac{2\res\varphi(k)}{u-k}\varphi(2u-k)-\frac{4\pi\gamma}{\lambda}\sum_{n=1}^{\infty}\sigma_{\alpha,k-2u}(n) n^{\frac{1-k}{2}}\int_{0}^{1}t^{\frac{k-1}{2}-u}J_{k-1}\left( \frac{4\pi\sqrt{nt}}{\lambda} \right)dt.
     \end{align*}
\end{proposition}

Here we trivially have 
\begin{align*}
    \res\varphi(k)=\left( \frac{2\pi}{\lambda} \right)^k \frac{\gamma\alpha_0}{\Gamma(k)}.
\end{align*}

\subsection{L-Functions Associated to Modular Forms}
Let $k\ge 4$ be an even integer and $\varphi$ be a Hecke series of signature $\left(1,k,(-1)^{k/2}\right)$. Then $\varphi$ is the L-function associated with $f_\alpha$ which is a modular form of weight $k$. We note that if $f_\alpha$ is a cusp form (if $\alpha_0=0$) then $c=k$ and $c=2k-1$ if $f_\alpha$ is not a cusp form \cite{ApostolModular}. Thus we have the following.
\begin{proposition}\label{CuspFormProp}
    If $f_\alpha$ is a cusp form, then for $\re(u),\re(v)>k+1$, we have
    \begin{align*}
        \varphi(u)\varphi(v)&=(-1)^{1+k/2}2\pi\sum_{n=1}^{\infty}\sigma_{\alpha,k-u-v}(n) n^{\frac{1-k}{2}}\int_{0}^{1}t^{\frac{k-1}{2}-u}J_{k-1}\left( 4\pi\sqrt{nt} \right)dt\\
        &+(-1)^{1+k/2}2\pi\sum_{n=1}^{\infty}\sigma_{\alpha,k-u-v}(n) n^{\frac{1-k}{2}}\int_{0}^{1}t^{\frac{k-1}{2}-v}J_{k-1}\left( 4\pi\sqrt{nt} \right)dt.
    \end{align*}
    Furthermore, for $\re(u)>k+1$, we have
    \begin{align*}
        \varphi(u)^2=(-1)^{1+k/2}4\pi\sum_{n=1}^{\infty}\sigma_{\alpha,k-2u}(n) n^{\frac{1-k}{2}}\int_{0}^{1}t^{\frac{k-1}{2}-u}J_{k-1}\left( 4\pi\sqrt{nt} \right)dt.
    \end{align*}
\end{proposition}
\begin{proposition}\label{NonCuspFormProp}
    If $f_\alpha$ is not a cusp form, then for $\re(u),\re(v)>2k$, we have
    \begin{align*}
        \varphi(u)\varphi(v)&=\frac{(2\pi i)^k\alpha_0}{(k-1)!}\left( \frac{1}{u-k}+\frac{1}{v-k} \right)\varphi(u+v-k)\\
        &+(-1)^{1+k/2}2\pi\sum_{n=1}^{\infty}\sigma_{\alpha,k-u-v}(n) n^{\frac{1-k}{2}}\int_{0}^{1}t^{\frac{k-1}{2}-u}J_{k-1}\left( 4\pi\sqrt{nt} \right)dt\\
        &+(-1)^{1+k/2}2\pi\sum_{n=1}^{\infty}\sigma_{\alpha,k-u-v}(n) n^{\frac{1-k}{2}}\int_{0}^{1}t^{\frac{k-1}{2}-v}J_{k-1}\left( 4\pi\sqrt{nt} \right)dt.
    \end{align*}
    Furthermore, for $\re(u)>2k$, we have
    \begin{align*}
        \varphi(u)^2=\frac{2(2\pi i)^k\alpha_0}{(k-1)!(u-k)}\varphi(2u-k)+(-1)^{1+k/2}4\pi\sum_{n=1}^{\infty}\sigma_{\alpha,k-2u}(n) n^{\frac{1-k}{2}}\int_{0}^{1}t^{\frac{k-1}{2}-u}J_{k-1}\left( 4\pi\sqrt{nt}\right)dt.
    \end{align*}
\end{proposition}
As an example of a cusp form we consider the modular discriminant $\Delta$, defined on $\tau\in\mathbb{H}$ by
\begin{align*}
    \Delta(\tau)=(2\pi)^{12}e^{2\pi i\tau}\prod_{n=1}^{\infty}\left(1-e^{2\pi i n\tau}\right)^{24}.
\end{align*}
$\Delta$ is a modular form of weight $12$ and has the Fourier expansion
\begin{align*}
    \Delta(\tau)=(2\pi)^{12}\sum_{n=1}^{\infty}\tau(n)e^{2\pi in \tau},
\end{align*}
where $\tau(n)$ is the Ramanujan's $\tau$ function. The associated Ramanujan's L-function \cite{BerndtHeckeTheory} is defined as
\begin{align*}
    L_\tau(s)=\sum_{n=1}^{\infty}\frac{\tau(n)}{n^s} \qquad \left( \sigma>13/2 \right).
\end{align*}
$L_\tau$ has the Hecke signature $(1,12,1)$ and can be analytically continued to an entire function. We get the following corollary of the Proposition \ref{CuspFormProp}.
\begin{corollary}
    For $\re(u),\re(v)>13$, we have
    \begin{align*}
        L_\tau(u)L_\tau(v)=&-2\pi\sum_{n=1}^{\infty}\sigma_{\tau,12-u-v}(n) n^{-11/2}\int_{0}^{1}t^{11/2-u}J_{11}\left( 4\pi\sqrt{nt}\right)dt\\
        &-2\pi\sum_{n=1}^{\infty}\sigma_{\tau,12-u-v}(n) n^{-11/2}\int_{0}^{1}t^{11/2-v}J_{11}\left( 4\pi\sqrt{nt} \right)dt.
    \end{align*}
    Furthermore, for $\re(u)>13$, we have
    \begin{align*}
        L_\tau(u)^2=-4\pi\sum_{n=1}^{\infty}\sigma_{\tau,12-2u}(n) n^{-11/2}\int_{0}^{1}t^{11/2-u}J_{11}\left( 4\pi\sqrt{nt} \right)dt.
    \end{align*}
\end{corollary}
As an example of a non-cusp form, we give the normalized Eisenstein series of order $k$ defined on $\tau\in\mathbb{H}$ by
\begin{align*}
    E_k(\tau)=1-\frac{2k}{B_k}\sum_{n=1}^{\infty}\sigma_{k-1}(n)e^{2\pi in\tau},
\end{align*}
where $k\ge4$ is an even integer and $B_k$ is the $k$th Bernoulli number. $E_k$ is a modular form of weight $k$. The associated L-function is given by
\begin{align*}
    L_k(s)=-\frac{2k}{B_k}\sum_{n=1}^{\infty}\frac{\sigma_{k-1}(n)}{n^s}=-\frac{2k}{B_k}\zeta(s)\zeta(s-k+1)\qquad \left( \sigma>k \right)
\end{align*}
and has the Hecke signature $\left(1,k,(-1)^{k/2}\right)$. Using the well-known formula $B_k=-\frac{2k!\zeta(k)}{(2\pi i)^k}$, we obtain the following $4$ term product identity by Proposition \ref{NonCuspFormProp}.
\begin{corollary}
    For $\re(u),\re(v)>2k$, we have
    \begin{align*}
        \zeta(u)\zeta(v)\zeta(u-k+1)\zeta(v-k+1)&=\zeta(k)\left( \frac{1}{u-k}+\frac{1}{v-k} \right)\zeta(u+v-k)\zeta(u+v-2k+1)\\
        &+(-1)^{1+k/2}2\pi\sum_{n=1}^{\infty}\sigma_{k,k-u-v}^*(n) n^{\frac{1-k}{2}}\int_{0}^{1}t^{\frac{k-1}{2}-u}J_{k-1}\left(4\pi\sqrt{nt} \right)dt\\
        &+(-1)^{1+k/2}2\pi\sum_{n=1}^{\infty}\sigma_{k,k-u-v}^*(n) n^{\frac{1-k}{2}}\int_{0}^{1}t^{\frac{k-1}{2}-v}J_{k-1}\left( 4\pi\sqrt{nt} \right)dt
    \end{align*}
    where $\sigma_{k,z}^*(n)=\sum_{d|n}\sigma_{k-1}(n)\sigma_{k-1}(n/d)d^z$. Furthermore, for $\re(u)>2k$, we have
    \begin{align*}
        \zeta(u)^2\zeta(u-k+1)^2&=\frac{2\zeta(k)}{u-k}\zeta(2u-k)\zeta(2u-2k+1)\\
        &+(-1)^{1+k/2}4\pi\sum_{n=1}^{\infty}\sigma_{k,k-2u}^*(n) n^{\frac{1-k}{2}}\int_{0}^{1}t^{\frac{k-1}{2}-u}J_{k-1}\left( 4\pi\sqrt{nt} \right)dt.
    \end{align*}
\end{corollary}
\subsection{Epstein Zeta Functions}
Let $Q$ be a positive definite quadratic form of $m$ variables, real coefficients, discriminant $D>0$ and inverse $Q^{-1}$. Then the Epstein zeta function \cite{BerndtHeckeTheory} is defined as
\begin{align*}
    Z(s;Q)=\sum_{\substack{\mathbf{n}\in \mathbb{Z}^m\\ \mathbf{n}\neq0}}\frac{1}{Q(\mathbf{n})^s}, \qquad \left( \sigma>m/2 \right).
\end{align*}
It satisfies the functional equation \cite{Epstein}
\begin{align*}
    \pi^{-s}\Gamma(s)Z(s;Q)=D^{-1/2}\pi^{s-m/2}\Gamma\left( \frac{m}{2}-s \right)Z\left( \frac{m}{2} -s;Q^{-1}\right).
\end{align*}
Thus it may be analytically continued to whole complex plane except $s=m/2$ where it has a pole with the residue
\begin{align*}
    \res Z\left(\frac{m}{2};Q\right)= D^{-1/2}\frac{\pi^{m/2}}{\Gamma(m/2)}.
\end{align*}
Assume that $Q$ has integer coefficients, then we get
\begin{align*}
    Z(s;Q)=\sum_{n=1}^{\infty}\frac{r_Q(n)}{n^s},
\end{align*}
where $r_Q(n)$ is the number of solutions to the equation $Q(\mathbf{n})=n$. Thus $\varphi(s)=Z(s;Q)$, $\psi(s)=Z(s;Q^{-1})$ and $\left(\lambda,k,\gamma\right)=\left(2,m/2,D^{-1/2}\right)$ satisfy the conditions of Theorem \ref{HCT}. Therefore we have the following.
\begin{corollary}\label{EpsteinCor}
    For $\re(u),\re(v)>1+m/2$, we have
    \begin{align*}
        Z(u;Q)Z(v;Q^{-1})&=\frac{\pi^{m/2}}{D^{1/2}\Gamma(m/2)(u-k)}Z(u+v-k;Q^{-1})+\frac{\pi^{m/2}D^{1/2}}{\Gamma(m/2)(v-k)}Z(u+v-k;Q)\\
        &-\frac{\pi}{D^{1/2}}\sum_{n=1}^{\infty}\sigma_{Q^{-1},m/2-u-v}(n) n^{\frac{2-m}{4}}\int_{0}^{1}t^{\frac{m-2}{4}-u}J_{m/2-1}\left( 2\pi\sqrt{nt} \right)dt\\
        &-\pi D^{1/2}\sum_{n=1}^{\infty}\sigma_{Q,m/2-u-v}(n) n^{\frac{2-m}{4}}\int_{0}^{1}t^{\frac{m-2}{4}-v}J_{m/2-1}\left( 2\pi\sqrt{nt}\right)dt,
    \end{align*}
    where $\sigma_{Q,z}(n)=\sigma_{r_{Q},z}(n)$. Furthermore, for $\re(u)>1+m/2$, we have
    \begin{align*}
        Z(u;Q)Z(u;Q^{-1})&=\frac{\pi^{m/2}}{D^{1/2}\Gamma(m/2)(u-k)}Z(2u-k;Q^{-1})+\frac{\pi^{m/2}D^{1/2}}{\Gamma(m/2)(u-k)}Z(2u-k;Q)\\
        &-\frac{\pi}{D^{1/2}}\sum_{n=1}^{\infty}\sigma_{Q^{-1},m/2-2u}(n) n^{\frac{2-m}{4}}\int_{0}^{1}t^{\frac{m-2}{4}-u}J_{m/2-1}\left( 2\pi\sqrt{nt} \right)dt\\
        &-\pi D^{1/2}\sum_{n=1}^{\infty}\sigma_{Q,m/2-2u}(n) n^{\frac{2-m}{4}}\int_{0}^{1}t^{\frac{m-2}{4}-u}J_{m/2-1}\left( 2\pi\sqrt{nt}\right)dt,
    \end{align*}
\end{corollary}
We give special attention to the case where $Q=Q^{-1}$. Then $Z(s;Q)$ is a Hecke series of signature $(2,m/2,1)$ and we get the following.
\begin{corollary}\label{EpsteinCor2}
    If $Q=Q^{-1}$, for $\re(u),\re(v)>1+m/2$, we have
    \begin{align*}
        Z(u;Q)Z(v;Q)&=\frac{\pi^{m/2}}{\Gamma(m/2)}\left( \frac{1}{u-k}+\frac{1}{v-k} \right)Z(u+v-k;Q)\\
        &-\pi\sum_{n=1}^{\infty}\sigma_{Q,m/2-u-v}(n) n^{\frac{2-m}{4}}\int_{0}^{1}t^{\frac{m-2}{4}-u}J_{m/2-1}\left( 2\pi\sqrt{nt} \right)dt\\
        &-\pi \sum_{n=1}^{\infty}\sigma_{Q,m/2-u-v}(n) n^{\frac{2-m}{4}}\int_{0}^{1}t^{\frac{m-2}{4}-v}J_{m/2-1}\left( 2\pi\sqrt{nt}\right)dt,
    \end{align*}
    Furthermore, for $\re(u)>1+m/2$, we have
    \begin{align*}
        Z(u;Q)^2=\frac{2\pi^{m/2}}{\Gamma(m/2)(u-k)}Z(2u-k;Q)-2\pi\sum_{n=1}^{\infty}\sigma_{Q,m/2-2u}(n) n^{\frac{2-m}{4}}\int_{0}^{1}t^{\frac{m-2}{4}-u}J_{m/2-1}\left( 2\pi\sqrt{nt} \right)dt
    \end{align*}
\end{corollary}

We note that the positive definite quadratic form $Q_m(n_1,\ldots,n_m)=n_1^2+\ldots+n_m^2$ satisfies the conditions of Corollary \ref{EpsteinCor2}.
\subsection{Dedekind Zeta Functions}
Wilton's formula for Dedekind zeta functions of both real and imaginary quadratic rings was investigated in \cite{DedekindWilton2,DedekindWilton3,DedekindWilton} and was applied to calculation of special values of Dedekind zeta functions. Unfortunately, Dedekind zeta function is a Hecke series only when $K$ is an imaginary quadratic field. We see that special cases of Theorem \ref{MainThm} indeed reduce to the results in literature.
\newline

Let $K$ be an imaginary quadratic field with the ring of integers $\mathcal{O}_K$, discriminant $d_K$, class number $h_K$ and number of roots of unity in $K$ $w_K$. Let $N(\mathfrak{a})=\left[ \mathcal{O}_K:\mathfrak{a} \right]$ denote the norm of an ideal $\mathfrak{a}$ in $\mathcal{O}_K$. We define the Dedekind zeta function \cite{DedekindWilton} by
\begin{align*}
    \zeta_K(s)=\sum_{\mathfrak{a}}\frac{1}{N(\mathfrak{a})^s}=\sum_{n=1}^{\infty}\frac{v_K(n)}{n^s}\qquad (\sigma>1),
\end{align*}
where the summation runs through the non-zero ideals of $\mathcal{O}_K$ and $v_K(n)$ denotes the number of ideals in $\mathcal{O}_K$ with norm $n$. Then by the well-known functional equation of the Dedekind zeta function, $\zeta_K$ is a Hecke series of signature $\left(|d_K|^{1/2},1,1\right)$ \cite{BerndtHeckeTheory}. Regulator of an imaginary quadratic field is equal to $1$. Therefore, the analytic class number formula reads
\begin{align*}
    \res\zeta_K(1)=\frac{2\pi h_K}{w_K |d_K|^{1/2}}.
\end{align*}
Thus, we have the following.
\begin{corollary}\label{DedekindCor}
    For $\re(u),\re(v)>1$, and $u,v\neq 2$, we have
    \begin{align*}
        \zeta_K(u)\zeta_K(v)&=\frac{2\pi h_K}{w_K\left| d_K \right|^{1/2}}\left( \frac{1}{u-k}+\frac{1}{v-k} \right)\zeta_K(u+v-1)\\
        &-\frac{2\pi}{|d_K|^{1/2}}\sum_{n=1}^{\infty}\sigma_{K,1-u-v}(n) \int_{0}^{1}t^{-u}J_{0}\left( 4\pi\sqrt{\frac{{nt}}{|d_k|}} \right)dt\\
        &-\frac{2\pi}{|d_K|^{1/2}}\sum_{n=1}^{\infty}\sigma_{K,1-u-v}(n) \int_{0}^{1}t^{-v}J_{0}\left( 4\pi\sqrt{\frac{{nt}}{|d_k|}} \right)dt,
    \end{align*}
    where $\sigma_{K,z}(n)=\sigma_{v_K,z}(n)$. Furthermore, for $\re(u)>1$, and $u\neq 2$, we have
    \begin{align*}
        \zeta_K(u)^2=\frac{4\pi h_K}{w_K\left| d_K \right|^{1/2}(u-k)}\zeta_K(2u-1)-\frac{4\pi}{|d_K|^{1/2}}\sum_{n=1}^{\infty}\sigma_{K,1-2u}(n) \int_{0}^{1}t^{-u}J_{0}\left( 4\pi\sqrt{\frac{{nt}}{|d_k|}} \right)dt.
    \end{align*}
\end{corollary}
Corollary \ref{DedekindCor} is just a mere rearrangement of the results given in \cite{DedekindWilton2,DedekindWilton3,DedekindWilton}, with less optimal conditions on $u,v$.

\subsection{Dirichlet L-Functions}
Let $\chi$ be a  primitive Dirichlet character $\bmod q$ where $q>1$. Let the Dirichlet L-function of $\chi$ be defined as
\begin{align*}
    L(s,\chi)=\sum_{n=1}^{\infty} \frac{\chi(n)}{n^s}\qquad (\sigma>1).
\end{align*}
First assume that $\chi$ is an even character. The functional equation for $L(s,\chi)$ \cite{Davenport} is given as follows,
\begin{align*}
    \pi^{-s/2}q^{s/2}\Gamma\left( \frac{s}{2} \right)L(s,\chi)=\frac{\tau(\chi)}{q^{1/2}}\pi^{(s-1)/2}q^{(1-s)/2}\Gamma\left( \frac{1-s}{2} \right)L(1-s,\overline{\chi}),
\end{align*}
where $\tau(\chi)=\sum_{n=1}^q \chi(n)\exp(2\pi i n/q)$ is the Gauss sum associated to $\chi$. This shows that $\varphi(s)=L(2s,\chi)$, $\psi(s)=L(2s,\overline{\chi})$ and $\left(\lambda,k,\gamma\right)=\left(2q,\frac{1}{2},\frac{\tau(\chi)}{q^{1/2}}\right)$ satisfy the conditions of Theorem \ref{HCT}. Since $L$-function of a non-principal character is entire, we get the following.
\begin{corollary}
    If $\chi$ is an even character, for $\re(u),\re(v)>1$, and $u,v\neq 3$, we have
    \begin{align*}
        L(u,\chi)L(v,\overline\chi)=&-\frac{\pi\tau(\chi)}{q^{3/2}}\sum_{n=1}^{\infty}\overline{\chi(n)}\sigma_{1-u-v}(n) n^{1/2}\int_{0}^{1}t^{-1/4-u/2}J_{-1/2}\left( \frac{2\pi n\sqrt{t}}{q} \right)dt\\
        &-\frac{\pi}{q^{1/2}\tau(\chi)}\sum_{n=1}^{\infty}\chi(n)\sigma_{1-u-v}(n) n^{1/2}\int_{0}^{1}t^{-1/4-v/2}J_{-1/2}\left( \frac{2\pi n\sqrt{t}}{q} \right)dt.
    \end{align*}
    Furthermore, for $\re(u)>1$ and $u\neq 3$, we have
    \begin{align*}
        L(u,\chi)L(u,\overline\chi)=&-\frac{\pi\tau(\chi)}{q^{3/2}}\sum_{n=1}^{\infty}\overline{\chi(n)}\sigma_{1-2u}(n) n^{1/2}\int_{0}^{1}t^{-1/4-u/2}J_{-1/2}\left( \frac{2\pi n\sqrt{t}}{q} \right)dt\\
        &-\frac{\pi}{q^{1/2}\tau(\chi)}\sum_{n=1}^{\infty}\chi(n)\sigma_{1-2u}(n) n^{1/2}\int_{0}^{1}t^{-1/4-u/2}J_{-1/2}\left( \frac{2\pi n\sqrt{t}}{q} \right)dt.
    \end{align*}
\end{corollary}
We note that the above formula is particularly useful when $\chi$ is a real character and $\varphi(s)=L(2s,\chi)$ is a Hecke series of signature $\left(2q,\frac{1}{2},\frac{\tau(\chi)}{q^{1/2}}\right)$ as a consequence. Furthermore, the Bessel function of order $-1/2$ may be computed explicitly as
\begin{align*}
    J_{-1/2}(z)=\sqrt{\frac{2}{\pi z}}\cos z.
\end{align*}

Now let $\chi$ be an odd character. Similarly, the functional equation for $L(s,\chi)$ \cite{Davenport} is given by
\begin{align*}
        \pi^{-(s+1)/2}q^{(s+1)/2}\Gamma\left( \frac{s+1}{2} \right)L(s,\chi)=\frac{\tau(\chi)}{iq^{1/2}}\pi^{(s-2)/2}q^{(2-s)/2}\Gamma\left( \frac{2-s}{2} \right)L(1-s,\overline{\chi}).
\end{align*}
This shows that $\varphi(s)=L(2s-1,\chi)$, $\psi(s)=L(2s-1,\overline{\chi})$ and $\left(\lambda,k,\gamma\right)=\left(2q,\frac{3}{2},\frac{\tau(\chi)}{iq^{1/2}}\right)$ satisfy the conditions of Theorem \ref{HCT}.
\begin{corollary}
    If $\chi$ is an odd character, for $\re(u),\re(v)>2$,  we have
    \begin{align*}
        L(u,\chi)L(v,\overline\chi)&=\frac{i\pi\tau(\chi)}{q^{3/2}}\sum_{n=1}^{\infty}\overline{\chi(n)}\sigma_{1-u-v}(n) n^{-3/2}\int_{0}^{1}t^{-1/4-u/2}J_{1/2}\left( \frac{2\pi n\sqrt{t}}{q} \right)dt\\
        &-\frac{i\pi}{q^{1/2}\tau(\chi)}\sum_{n=1}^{\infty}\chi(n)\sigma_{1-u-v}(n) n^{-3/2}\int_{0}^{1}t^{-1/4-v/2}J_{1/2}\left( \frac{2\pi n\sqrt{t}}{q} \right)dt.
    \end{align*}
    Furthermore, for $\re(u)>2$, we have
    \begin{align*}
        L(u,\chi)L(u,\overline\chi)&=\frac{i\pi\tau(\chi)}{q^{3/2}}\sum_{n=1}^{\infty}\overline{\chi(n)}\sigma_{1-2u}(n) n^{-3/2}\int_{0}^{1}t^{-1/4-u/2}J_{1/2}\left( \frac{2\pi n\sqrt{t}}{q} \right)dt\\
        &-\frac{i\pi}{q^{1/2}\tau(\chi)}\sum_{n=1}^{\infty}\chi(n)\sigma_{1-2u}(n) n^{-3/2}\int_{0}^{1}t^{-1/4-u/2}J_{1/2}\left( \frac{2\pi n\sqrt{t}}{q} \right)dt.
    \end{align*}
\end{corollary}
Here again the above formula is particularly useful when $\chi$ is a real character and $\varphi(s)=L(2s-1,\chi)$ is a Hecke series of signature $\left(2q,\frac{3}{2},\frac{\tau(\chi)}{i q^{1/2}}\right)$ as a consequence. Similarly, the Bessel function of order $1/2$ may be computed explicitly as
\begin{align*}
    J_{1/2}(z)=\sqrt{\frac{2}{\pi z}}\sin z.
\end{align*}

\subsection{Riemann Zeta Function}
By the equation \eqref{RiemannFe}, it is evident that $\varphi(s)=\zeta(2s)$ is a Hecke series of signature $\left(2,1/2,1\right)$. We now give a slightly different version of Wilton's formula.
\begin{corollary}
    For $\re(u),\re(v)>1$, and $u,v\neq 3$, we have
    \begin{align*}
        \zeta(u)\zeta(v)-\left( \frac{1}{u-1}+\frac{1}{v-1} \right)\zeta(u+v-1)=
        &-\pi\sum_{n=1}^{\infty}\sigma_{1-u-v}(n) n^{1/2}\int_{0}^{1}t^{-1/4-u/2}J_{-1/2}\left(2\pi n\sqrt{t} \right)dt\\
        &-\pi\sum_{n=1}^{\infty}\sigma_{1-u-v}(n) n^{1/2}\int_{0}^{1}t^{-1/4-v/2}J_{-1/2}\left(2\pi n\sqrt{t} \right)dt.
    \end{align*}
    Furthermore, for $\re(u)>1$ and $u\neq3$, we have
    \begin{align*}
        \zeta(u)^2-\frac{2}{u-1}\zeta(2u-1)=-2\pi\sum_{n=1}^{\infty}\sigma_{1-2u}(n) n^{1/2}\int_{0}^{1}t^{-1/4-u/2}J_{-1/2}\left(2\pi n\sqrt{t} \right)dt.
    \end{align*}
\end{corollary}

\bibliography{main}
\end{document}